\documentclass[11pt]{article}
\usepackage[T1]{fontenc}
\usepackage{lmodern,microtype}
\usepackage[margin=1in]{geometry}
\usepackage{amsmath,amssymb,amsthm,mathtools,enumitem,booktabs,xcolor}
\usepackage[colorlinks=true,linkcolor=blue!40!black,citecolor=blue!40!black,urlcolor=blue!40!black]{hyperref}
\usepackage{xurl}
\hypersetup{pdftitle={Deterministic Polynomial-Time Primality Testing with a Single Binomial Congruence},pdfauthor={Joseph M. Shunia},pdfsubject={Deterministic primality testing in cyclotomic and binomial quotient algebras}}
\newtheorem{theorem}{Theorem}[section]
\newtheorem{lemma}[theorem]{Lemma}
\newtheorem{proposition}[theorem]{Proposition}
\newtheorem{corollary}[theorem]{Corollary}
\newtheorem{conjecture}[theorem]{Conjecture}
\theoremstyle{definition}

\theoremstyle{remark}
\newtheorem{remark}[theorem]{Remark}
\newcommand{\Z}{\mathbb Z}
\newcommand{\Q}{\mathbb Q}
\newcommand{\F}{\mathbb F}
\newcommand{\ord}{\operatorname{ord}}
\newcommand{\softO}{\widetilde O}
\newcommand{\A}{\mathcal A}
\setlist[enumerate]{itemsep=4pt,topsep=5pt}
\title{Deterministic Polynomial-Time Primality Testing\\with a Single Binomial Congruence}
\author{Joseph M. Shunia}
\date{October 30, 2026}

\begin{document}
\maketitle

\begin{abstract}
We give a deterministic polynomial-time primality test whose main
verification is the congruence $(1+X)^n=1+X^n$ in $C_n[X]/(X^r-b(Y))$, where
$C_n=(\Z/n\Z)[Y]/(\Phi_r(Y))$. Two identities in the coefficient
ring normalize $b$ and ensure compatibility with the main congruence.
A finite-field root count and a lattice bound show that an accepted
input is a prime power; elementary preprocessing excludes nontrivial
perfect powers. We construct $b$ deterministically using Jacobi sums
and the uniform Dedekind-zeta zero-free theorem of Cureton and Thorner,
which extends OpenAI's Hecke zero-free theorem. The total bit complexity
is $L^{32+o(1)}$, where $L=\lceil\log_2 n\rceil$, and verification
after constructing $b$ costs $L^{6+o(1)}$. More generally, for a fixed
real algebraic $1/2\le\theta<1$, a common zero-free boundary
$\Re s>\theta$ for the relevant Dedekind zeta functions gives
$L^{\max\{6,(1-\theta)^{-1}\}+o(1)}$. This bound reaches
$L^{6+o(1)}$ at $\theta=5/6$, and hence under GRH. We state a fixed-base modification of
Agrawal's conjecture for odd integers passing the base-$2$ Fermat
test, using the modulus $X^r-2$ with $r$ the least positive integer
not dividing $n(n-1)$.
\end{abstract}

\section{Introduction}

For a prime $n$, the identity
\begin{equation}\label{eq:frobenius-intro}
 (1+X)^n=1+X^n
\end{equation}
holds in every commutative algebra of characteristic $n$. A primality
test based on this identity must choose an algebra small enough to
compute in and sufficiently restrictive to exclude composite inputs.
We construct such an algebra deterministically, using one shift, one
binomial relation, and exponent $n$.

The coefficient ring and the final algebra are
\begin{equation}\label{eq:rings-intro}
 C_n=(\Z/n\Z)[Y]/(\Phi_r(Y)),\qquad
 \A_n=C_n[X]/(X^r-b(Y)),
\end{equation}
where $r$ is odd and squarefree and $\Phi_r$ is the $r$th cyclotomic
polynomial. Thus the test can equally be written in the multivariate
quotient
\[
 (\Z/n\Z)[X,Y]/\bigl(\Phi_r(Y),X^r-b(Y)\bigr).
\]
Arithmetic is performed in the full quotient, which may have zero
divisors. Field quotients enter only in the proof.

For a suitable auxiliary integer $D$, the acceptance equalities are
\begin{align}
 b(Y)^{(n^D-1)/r}&=Y &&\text{in } C_n,\label{eq:normalize}\\
 b(Y)^n&=b(Y^n) &&\text{in } C_n,\label{eq:compat}\\
 (1+X)^n&=1+X^n &&\text{in }\A_n.\label{eq:main}
\end{align}
The single binomial congruence is \eqref{eq:main}, evaluated in
$\A_n$. The two checks in $C_n$ and elementary integer preprocessing
complete the test.

The deterministic parameter construction uses a uniform zero-free
theorem. OpenAI proved a common width $10^{-6}$ for finite-order
Hecke $L$-functions over cyclotomic fields containing $\mu_{12}$
\cite[Theorem 1.2]{openai-artin}. Cureton and Thorner extended this
approach to all number fields and obtained the following stronger
result.

\begin{theorem}[Cureton--Thorner {\cite[Theorem 1.1]{cureton-thorner}}]
\label{thm:zero-free}
For every number field $F$, its Dedekind zeta function satisfies
\[
 \zeta_F(s)\ne0\qquad
 \left(\Re s>\frac{31}{32},\ s\ne1\right).
\]
\end{theorem}

This theorem has a published Lean formalization generated with
AxiomProver \cite{axiom-lean}. We use $\delta=1/32$ in the main
algorithm, so that its search cutoff is computed by integer powering.
We retain $\delta$ in the estimates to exhibit the effect of stronger
analytic bounds. The uniformity needed below concerns the zero-free
boundary; Lemma~\ref{lem:explicit} separately makes the effective
dependence on field degree and discriminant explicit.

The main result uses Theorem~\ref{thm:zero-free} as its analytic
input. Corollary~\ref{cor:analytic-complexity} gives the dependence
of the running time on a general uniform zero-free boundary.
Any stronger boundary satisfying its hypotheses therefore gives
the corresponding improved bound by direct substitution, with
the same parameter construction and algebraic acceptance criterion.
Section~\ref{sec:conditional-bounds} also identifies the $5/6$
threshold and the GRH endpoint.

\begin{theorem}[Main result]\label{thm:main}
There is a deterministic primality algorithm
using \eqref{eq:normalize}--\eqref{eq:main} with total bit complexity
\[
 L^{32+o(1)},\qquad L=\lceil\log_2 n\rceil.
\]
Its parameters satisfy
\[
 D=L^{o(1)},\qquad r=L^{2+o(1)},\qquad
 \varphi(r)\le r.
\]
Verification after constructing $b$ costs $L^{6+o(1)}$ bit operations.
\end{theorem}

The proof has three components. A normalized coefficient first produces
a cyclic orbit of $r$ linear factors in a field extension. Bernstein's
root-counting and lattice method \cite{bernstein} then excludes inputs
with two distinct prime factors. Next, an Odlyzko--Pomerance construction
gives $r=L^{2+o(1)}$, with every prime factor of $r$ of size $L^{o(1)}$.
Finally, we construct the coefficient by combining Jacobi sums for those
small prime factors. Treating the prime degrees separately keeps the
number-field degrees in the auxiliary-prime estimate small.

The AKS theorem established deterministic polynomial-time primality
testing \cite{aks}. Berrizbeitia and Cheng developed tests using
binomial extensions and a small number of congruences
\cite{berrizbeitia,cheng}. Bernstein gave an essentially quartic
random-time primality-proving algorithm \cite{bernstein}; its
group-size argument is the algebraic foundation used here.
Lenstra and Pomerance subsequently obtained an unconditional
$L^{6+o(1)}$ deterministic algorithm using Gaussian periods
\cite{lenstra-pomerance}. The present construction concerns the
deterministic realization of the particular congruence
\eqref{eq:main}, with a normalized coefficient in the full cyclotomic
quotient. Its auxiliary-prime search dominates the stated total bound.

Throughout, $\log$ is the natural logarithm unless a base is displayed.
The notation $\softO$ suppresses factors polynomial in logarithms.
All complexity statements count bit operations; $o(1)$ is as
$L\to\infty$, with $\delta$ fixed.

\section{An algebraic acceptance criterion}\label{sec:criterion}

We first prove that the three acceptance identities force $n$ to be
a prime power. This part of the argument depends only on the identities
and the sizes of the degrees, not on how the coefficient is constructed.

Let $n>3$ be odd, and choose an odd squarefree integer $r>2$ and a
positive integer $D$ such that
\begin{equation}\label{eq:degrees}
 \gcd(n,r)=\gcd(D,r)=1,\qquad q-1\mid D
       \quad\text{for every prime }q\mid r.
\end{equation}
In particular, $\lambda(r)\mid D$, so $n^D\equiv1\pmod r$ and the
exponent in \eqref{eq:normalize} is integral. For odd $r$, define
\begin{equation}\label{eq:B}
 h=\frac{r-1}{2},\qquad
 B_r=\sum_{j=0}^{h}\binom rj\binom hj\binom{r-j+h}{h}.
\end{equation}
Our exact size condition is
\begin{equation}\label{eq:size}
 3\bigl\lfloor\log_2 B_r\bigr\rfloor^2\ge D^2L^2r,
 \qquad L=\lceil\log_2 n\rceil.
\end{equation}
This condition uses integer arithmetic only: the logarithm floor is
the bit length of $B_r$ minus one.

\begin{lemma}[Cyclotomic coefficient ring]\label{lem:ring}
Under \eqref{eq:degrees}, substitution $Y\mapsto Y^n$ defines an
automorphism $\sigma$ of $C_n$ with $\sigma^D=1$. Each field quotient
of $C_n$ is $\F_{p^e}$ for some prime $p\mid n$, where
$e=\ord_r(p)\mid D$. In every such quotient the image of $Y$ has
order exactly $r$.
\end{lemma}

\begin{proof}
Since $\gcd(n,r)=1$, raising a primitive $r$th root to its $n$th power
permutes the primitive roots. The integral cyclotomic relation therefore
makes the substitution well-defined. An inverse is substitution by
$Y^v$ for $nv\equiv1\pmod r$. Also $Y^r=1$ and $n^D\equiv1\pmod r$,
which prove $\sigma^D=1$.

A maximal ideal of $C_n$ lies over a prime divisor $p$ of $n$.
Because $p\nmid r$, the polynomial $X^r-1$ is squarefree in
$\F_p[X]$. The roots of $\Phi_r$ are precisely the primitive
$r$th roots, whose Frobenius orbits have length $\ord_r(p)$.
The quotient fields consequently have the stated degrees. Finally,
$p^D\equiv1\pmod r$ implies $e\mid D$.
\end{proof}

\begin{lemma}[Iteration of the main congruence]\label{lem:iteration}
If \eqref{eq:compat} and \eqref{eq:main} hold, then
\[
 (1+X)^{n^j}=1+X^{n^j}\quad\text{in }\A_n
 \qquad(j\ge0).
\]
\end{lemma}

\begin{proof}
Equation \eqref{eq:compat} gives $\sigma^i(b)=b^{n^i}$ by induction.
Put $F(X)=(1+X)^n-1-X^n$. Equation \eqref{eq:main} says
$X^r-b\mid F(X)$ in $C_n[X]$. Applying $\sigma^i$ to its
coefficients gives $X^r-b^{n^i}\mid F(X)$. Substituting $X^{n^i}$
and using
\[
 X^r-b\mid X^{rn^i}-b^{n^i}
\]
gives $(1+X^{n^i})^n=1+X^{n^{i+1}}$ in $\A_n$. Starting at
$j=0$, this proves the assertion by induction. No additivity of
$n$th powering for composite $n$ is used.
\end{proof}

\begin{theorem}[Acceptance criterion]\label{thm:criterion}
Suppose \eqref{eq:degrees} and \eqref{eq:size} hold. If
\eqref{eq:normalize}--\eqref{eq:main} hold, then $n$ is a prime
power. If $n$ is not a nontrivial perfect power, it is prime.
\end{theorem}

\begin{proof}
Put $N=n^D$. The element $Y$ is a unit, so
\eqref{eq:normalize} makes $b$ a unit. Suppose for a contradiction
that $n$ has at least two distinct prime divisors. Choose a field
quotient $k=\F_{p^e}$ of $C_n$ as in Lemma~\ref{lem:ring}, and
put $P=p^e$. We use the same names for the images of $Y$ and $b$.
One has $P\mid N$, and $Y$ has order $r$ in $k^\times$.

Choose a root $\alpha$ of $X^r-b$ in an algebraic closure of $k$.
Normalization and Lemma~\ref{lem:iteration} imply
\[
 \alpha^N=Y\alpha,\qquad
 u_i:=1+Y^i\alpha,\qquad u_i^N=u_{i+1}\quad(i\bmod r).
\]
For the last identity, substitute $Y^i\alpha$ in
$(1+X)^N=1+X^N$ and use $N\equiv1\pmod r$.
Each $u_i$ is nonzero: otherwise $\alpha=-Y^{-i}$, so
$b=(-1)^r$ and $Y=b^{(N-1)/r}$ would have order at most two.
Let $G$ be the finite multiplicative group generated by the $u_i$
in $k(\alpha)^\times$.

We first prove $|G|\ge B_r$, using the signed-product count recorded
by Bernstein \cite[Section 3]{bernstein}. Consider all integer vectors
$a=(a_0,\ldots,a_{r-1})$ for which the positive entries sum to at
most $h$ and the absolute values of the negative entries also sum
to at most $h$. If exactly $j$ entries are negative, their positions
can be chosen in $\binom rj$ ways, their absolute values in
$\binom hj$ ways, and the remaining nonnegative entries in
$\binom{r-j+h}{h}$ ways. This gives exactly $B_r$ vectors.

Their products $\prod_i u_i^{a_i}$ are distinct. Indeed, suppose
two such products are equal. Taking successive $N$th powers rotates
all indices, so the equality holds with $\alpha$ replaced by each
of the $r$ distinct elements $Y^j\alpha$. Clearing denominators
produces an equality of two polynomials at all those elements;
each polynomial has degree at most $2h=r-1$. Their difference is
therefore zero by the field root bound. Unique factorization into
the pairwise distinct factors $1+Y^i Z$, whose constant terms are
all $1$, makes the original exponent vectors equal. Hence
\begin{equation}\label{eq:lower-group}
 |G|\ge B_r.
\end{equation}

For the upper bound, $r\mid P-1$, and $b\in k$ gives
$\alpha^P=Y^t\alpha$ for some integer $t$. Set $H=N/P>1$.
Frobenius in characteristic $p$, followed by injectivity of
$P$th powering, gives
\begin{equation}\label{eq:rotations}
 u_i^P=u_{i+t},\qquad u_i^H=u_{i+1-t}.
\end{equation}
Consider the lattice
\[
 \Lambda=\{(a,c)\in\Z^2:(1-t)a+tc\equiv0\pmod r\}.
\]
It has determinant $r$, because $1-t$ and $t$ generate the unit
ideal in $\Z/r\Z$. The symmetric convex hexagon
\[
 \max\{|a\log H|,|c\log P|,|a\log H+c\log P|\}
 \le \sqrt{r/3}\log N
\]
has area
\[
 \frac{r(\log N)^2}{\log H\log P}\ge4r,
\]
since $\log H+\log P=\log N$. Minkowski's convex body theorem
gives a nonzero lattice point in this closed hexagon. The equality
case follows by slightly enlarging the hexagon and using discreteness.
Change its sign if necessary to make $a\ge0$. If $c\ge0$, put
$U=H^aP^c$ and $V=1$; if $c<0$, put $U=H^a$ and $V=P^{-c}$.
Equation \eqref{eq:rotations} and the lattice relation give
$g^U=g^V$ for every $g\in G$, and
\[
 \max(U,V)\le N^{\sqrt{r/3}}.
\]
Moreover, $U\ne V$. For $c\ge0$ this follows from $(a,c)\ne(0,0)$.
For $c<0$, equality would contradict the fact that $H$ has a prime
factor other than $p$ when $a>0$; when $a=0$, it would give
$1=P^{-c}>1$.

Every element of $G$ is consequently a root of
$Z^{|U-V|}-1$ in a field. Thus
\[
 |G|\le |U-V|<N^{\sqrt{r/3}}\le B_r.
\]
The last inequality follows from \eqref{eq:size}, because
$\log_2 N\le DL$. This contradicts \eqref{eq:lower-group}.
Therefore $n$ has only one distinct prime divisor.
\end{proof}

\begin{remark}
The argument does not require $X^r-b$ to be irreducible modulo every
prime divisor of $n$. It takes a root in a field extension and obtains
the needed distinct products from the cyclic rotation and the field
root bound. Likewise, $C_n$ itself may have zero divisors. Its rank
$\varphi(r)$ measures arithmetic cost, whereas the smaller number
$e\mid D$ controls the lattice bound.
\end{remark}

\section{Choosing the degrees}\label{sec:degrees}

We first record the size of the signed-product count. The improvement
from signed products is due to Voloch and Vaaler, as explained in
Bernstein's analysis \cite[Section 3]{bernstein}.

\begin{proposition}\label{prop:count}
For odd $r\to\infty$,
\[
 \log_2 B_r=\kappa r+O(\log r),\qquad
 \kappa=\log_2(3+2\sqrt2).
\]
Consequently, for every fixed $\epsilon>0$,
\[
 r>(c_0+\epsilon)D^2L^2,\qquad
 c_0=\frac1{3\kappa^2}=0.051540616989\ldots,
\]
implies \eqref{eq:size} for sufficiently large $r$.
A simpler sufficient condition valid for odd $r\ge7$ is
$r\ge D^2L^2$.
\end{proposition}

\begin{proof}
In \eqref{eq:B}, retain the term with
$j=\lfloor(1-1/\sqrt2)r\rfloor$. Stirling's formula for its three
binomial coefficients gives $\log_2 B_r\ge\kappa r-O(\log r)$.
For the upper bound, set $t=\sqrt2-1$. The generating function for
one signed exponent, recording positive and negative weights
separately, is $1+z/(1-z)+w/(1-w)$. Every vector counted by $B_r$
has combined weight at most $2h$, so
\[
 B_r t^{2h}\le
 \left(1+\frac{2t}{1-t}\right)^r=t^{-r}.
\]
Since $2h=r-1$, this has the asserted exponential rate. The
asymptotic threshold follows on substituting this estimate in
\eqref{eq:size}. For the simpler bound, the count includes all
nonnegative $0$--$1$ vectors of weight at most $h$ and their
negatives, with zero counted once. Their number is $2^r-1$.
Hence $\lfloor\log_2 B_r\rfloor\ge r-1$, and
$3(r-1)^2\ge r^2\ge D^2L^2r$ for the stated range.
\end{proof}

To meet the size condition with small $D$, we build $r$ from primes
$q$ for which $q-1\mid D$. The following parameter theorem is due
to Odlyzko and Pomerance;
we use the form in Adleman--Pomerance--Rumely
\cite[Theorem 3]{apr}, also used by Bernstein \cite{bernstein}.
There is an absolute constant $C$ such that, for every sufficiently
large $H$, some positive integer $D$ satisfies
\begin{equation}\label{eq:OP}
 D\le\exp(C\log\log H\log\log\log H),\qquad
 \prod_{\substack{q-1\mid D\\q\text{ prime}}}q\ge H.
\end{equation}

\begin{proposition}\label{prop:degrees}
There is a deterministic search for $D$ and odd squarefree $r>2$
satisfying $\gcd(D,r)=1$, $q-1\mid D$ for every prime $q\mid r$,
and \eqref{eq:size}, with
\begin{equation}\label{eq:parameter-sizes}
 D=L^{o(1)},\qquad r=L^{2+o(1)},\qquad
 q\le D+1=L^{o(1)}\quad(q\mid r).
\end{equation}
Its cost is at most $L^{4+o(1)}$ bit operations.
\end{proposition}

\begin{proof}
Enumerate even positive integers $D$. For each, enumerate the primes
$q\le D+1$ with $q-1\mid D$ and $q\nmid D$. Multiply these primes
in increasing order, checking \eqref{eq:size} after each prefix
whose product exceeds $2$. If all primes are exhausted without
success, continue to the next $D$.

To bound the search, apply \eqref{eq:OP} with $H=L^3$. The supplied
$D$ is $L^{o(1)}$. Deleting the primes dividing $D$ removes a
product at most $D$, so the remaining product is at least
$H/D\ge D^2L^2$ for sufficiently large $L$. Such a successful
$D$ is even: if $D$ were odd, the product in \eqref{eq:OP} would
contain only the prime $2$. For even $D$, excluding $q\mid D$
also excludes $2$. Proposition~\ref{prop:count} now guarantees
success. For each fixed smaller $L$, applying \eqref{eq:OP} with
$H$ sufficiently large gives the same termination argument.

At the first successful prefix the product obeys
\[
 r\le(D+1)\max\{6,D^2L^2\}=L^{2+o(1)}.
\]
Conversely, Proposition~\ref{prop:count}'s upper estimate and
\eqref{eq:size} give $r\gg D^2L^2$, establishing the stated
asymptotic scale. Each trial has $L^{o(1)}$ small primes and
divisibility checks. The exact count $B_r$ can be evaluated from
binomial recurrences in $\softO(r^2)$ bit operations: there are
$O(r)$ terms, each with $O(r)$ bits. The $L^{o(1)}$ trials fit the
claimed bound.
\end{proof}

The factorization of $r$ is known from its construction. To enforce
$\gcd(n,r)=1$, test its known prime factors against $n$. A common
prime factor proves compositeness unless it equals $n$, in which
case $n$ is prime.

\section{A uniform auxiliary-prime bound}\label{sec:analytic}

For each prime $q\mid r$, we seek a prime $\ell$ modulo which $n$
is a $q$th power nonresidue. Such a prime will define the character
used in Section~\ref{sec:coefficient}. We prove its existence for
prime inputs by applying Theorem~\ref{thm:zero-free} to Kummer
extensions of degree $q$. The following smoothed estimate makes
the dependence on the number field explicit.

\begin{lemma}[Smoothed prime-ideal estimate]\label{lem:explicit}
Fix $0<\delta\le1/2$ and a nonzero nonnegative smooth function $w$
supported in $(1,2)$. Let $F$ be a number field of degree $d_F$ and absolute discriminant
$\Delta_F$, and suppose its Dedekind zeta function has no zeros in
$\Re s>1-\delta$. Put
\[
 \Psi_F(x)=\sum_{\mathfrak p}\sum_{a\ge1}
    \log N\mathfrak p\;w((N\mathfrak p)^a/x).
\]
Then, uniformly in $F$ and $x\ge2$,
\begin{equation}\label{eq:explicit}
 \Psi_F(x)=x\widehat w(1)
  +O_w\bigl(x^{1-\delta}(\log\Delta_F+d_F)\bigr),
 \qquad \widehat w(s)=\int_0^\infty w(u)u^{s-1}\,du.
\end{equation}
The implied constant is effective.
\end{lemma}

\begin{proof}
This is the smoothed explicit-formula calculation in
\cite[Section 9.2, especially (9.9)]{openai-artin}; we recall the
dependence on the field. Mellin inversion applied to $-\zeta_F'/\zeta_F$
and a contour shift from $\Re s=2$ to $\Re s=-1/2$ give
\[
 \Psi_F(x)=x\widehat w(1)
 -\sum_\rho x^\rho\widehat w(\rho)
 -(r_1+r_2-1)\widehat w(0)
 +O_w\bigl(x^{-1/2}(\log\Delta_F+d_F)\bigr),
\]
where $\rho$ runs over nontrivial zeros with multiplicity and
$(r_1,r_2)$ is the signature. The functional equation bounds the
integral on the left line with the displayed dependence. The
uniform zero count in a unit interval at height $T$ is
$O(\log\Delta_F+d_F\log(2+|T|))$. The Mellin transform decays
faster than every inverse power of $|T|$ in bounded vertical strips.
Since $\Re\rho\le1-\delta$, the zero sum is absolutely convergent
and is bounded by
$O_w(x^{1-\delta}(\log\Delta_F+d_F))$. The residue at zero and
the left-line integral are absorbed in the same bound. These
estimates have effective absolute constants for fixed $w$.
\end{proof}

\begin{lemma}[Auxiliary prime for one prime degree]\label{lem:aux}
Fix $0<\delta<1/2$, and suppose that for every number field $F$,
$\zeta_F(s)$ has no zeros in $\Re s>1-\delta$. There is an
effective constant
$C_\delta$ such that, for every prime $n>3$ and every odd prime
$q\ne n$, there is a prime $\ell$ satisfying
\begin{gather}
 \ell\le C_\delta\bigl(q(L+\log q+2)\bigr)^{1/\delta},
       \label{eq:aux-bound}\\
 \ell\equiv1\pmod{12q},\qquad \ell\nmid n,\qquad
 n^{(\ell-1)/q}\not\equiv1\pmod\ell.
       \label{eq:aux-properties}
\end{gather}
\end{lemma}

\begin{proof}
Set $K=\Q(\mu_{12q})$ and $M=K(n^{1/q})$. The prime $n$ is
unramified in $K$, so $X^q-n$ is Eisenstein at a prime of $K$
above $n$. Thus $[M:K]=q$, and $M/K$ is cyclic because
$\mu_q\subset K$. By hypothesis, both $\zeta_K$ and $\zeta_M$
are zero-free in $\Re s>1-\delta$. For the stated algorithm,
Theorem~\ref{thm:zero-free} supplies this hypothesis with
$\delta=1/32$.

The cyclotomic discriminant bound and the discriminant of $X^q-n$
give
\begin{align*}
 [K:\Q]&=O(q),& \log\Delta_K&=O(q\log q),\\
 [M:\Q]&=O(q^2),& \log\Delta_M&=O(q^2(L+\log q)).
\end{align*}
More explicitly, the relative discriminant divides
$(q^q n^{q-1})\mathcal O_K$, and the discriminant tower formula
gives the second bound. Lemma~\ref{lem:explicit} therefore yields
\begin{equation}\label{eq:difference}
 \Psi_K(x)-\frac1q\Psi_M(x)
 =\left(1-\frac1q\right)x\widehat w(1)
 +O_w\bigl(q(L+\log q+2)x^{1-\delta}\bigr).
\end{equation}

Consider an unramified degree-one prime of $K$. If it splits in
$M$, its first-power contribution cancels against the contributions
of the $q$ primes above it, after the factor $1/q$. Otherwise it
is inert, since the extension has prime degree, and its positive
first-power contribution remains. A surviving prime lies above a
rational prime $\ell\equiv1\pmod{12q}$ for which $n$ is not a
$q$th power in $\F_\ell^\times$. This is exactly the last
condition in \eqref{eq:aux-properties}.

All terms of prime-ideal residue degree at least two, and all
higher prime powers, contribute in absolute value at most
$O_w(q\sqrt{x}\log^2(2x))$ to \eqref{eq:difference}. To see the
field dependence, a field of degree $d$ contributes at most
$O_w(d\sqrt{x}\log^2(2x))$; the weighted total degree here is
$[K:\Q]+[M:\Q]/q=O(q)$. Ramified degree-one first powers
contribute at most $O_w(q\log(6qn))$, since the ramified rational
primes divide $6qn$.

Take
$x=T_\delta\bigl(q(L+\log q+2)\bigr)^{1/\delta}$.
For an effectively large constant $T_\delta$, the error in
\eqref{eq:difference} is less than half its main term. The
discarded terms are smaller as well, uniformly in $q$ and $n$,
because $\delta<1/2$. Some surviving prime lies in $(x,2x)$.
Absorbing $2T_\delta$ in $C_\delta$ proves the lemma.
\end{proof}

All primes $q\mid r$ are at most $D+1$. Hence a common search
bound, increased by rounding if needed, is
\begin{equation}\label{eq:search-bound}
 B=\left\lceil C_\delta
 \bigl((D+1)(L+\lceil\log_2(D+1)\rceil+3)\bigr)^{1/\delta}
 \right\rceil=L^{1/\delta+o(1)}=L^{32+o(1)}.
\end{equation}
For the search itself it suffices to test primes
$\ell\equiv1\pmod q$; Lemma~\ref{lem:aux} supplies a prime
even in the smaller progression. Only prime-input existence is
needed. If a search fails within the bound, the algorithm rejects.

\section{Constructing the normalized coefficient}\label{sec:coefficient}

\subsection{Characters defined by the input}

For each prime $q\mid r$, choose a prime $\ell_q\le B$ with
$q\mid\ell_q-1$, $\ell_q\nmid n$, and
$n^{(\ell_q-1)/q}\not\equiv1\pmod{\ell_q}$.
Put
\[
 t_q=n^{(\ell_q-1)/q}\bmod\ell_q.
\]
Its order in $\F_{\ell_q}^\times$ is exactly $q$. For $u\ne0$
define $j_q(u)\in\{0,\ldots,q-1\}$ by
\[
 u^{(\ell_q-1)/q}=t_q^{j_q(u)}\pmod{\ell_q}.
\]
This is a multiplicative-group homomorphism to $\Z/q\Z$, written
additively, and $j_q(n\bmod\ell_q)=1$. Its values are computed
by a table of the $q$ powers of $t_q$.

Choose $s_q\in\{1,\ldots,q-1\}$ satisfying
\begin{equation}\label{eq:sq}
 s_q\equiv-\bigl(D(r/q)\bigr)^{-1}\pmod q.
\end{equation}
The inverse exists by \eqref{eq:degrees} and squarefreeness. In
$R_q=(\Z/n\Z)[Z]/(\Phi_q(Z))$, define the order-$q$ character
\begin{equation}\label{eq:chi}
 \chi_q(u)=Z^{s_qj_q(u)}\quad(u\ne0),\qquad \chi_q(0)=0.
\end{equation}
In particular, $\chi_q(n)=Z^{s_q}$ and $\chi_q(-1)=1$, since
$q$ is odd. This construction needs neither a primitive root
modulo $\ell_q$ nor a factorization of $\ell_q-1$.

\subsection{Jacobi sums and a binary recurrence}

For positive integers $a,c$, set
\[
 J(\chi^a,\chi^c)=
 \sum_{u\in\F_\ell}\chi(u)^a\chi(1-u)^c.
\]
Here $\chi=\chi_q$ and $\ell=\ell_q$. The terms for $u=0,1$
vanish. The other terms are powers of $Z$, so the sum is computed
by counting $q$ possible exponents before reduction modulo
$(n,\Phi_q(Z))$. Define
\begin{equation}\label{eq:beta}
 \beta_q(Z)=\ell_q\prod_{j=1}^{q-2}J(\chi_q,\chi_q^j)
       \quad\text{in }R_q.
\end{equation}

\begin{lemma}\label{lem:jacobi}
The element $\beta_q$ is the $q$th power of its Gauss sum in a
cyclotomic extension of $R_q$. It can be computed in $R_q$ using
only $O(\log q)$ Jacobi sums, without constructing that extension.
\end{lemma}

\begin{proof}
For this proof only, work in $R_q[T]/(\Phi_\ell(T))$ and put
$\tau_a=\sum_{u\ne0}\chi(u)^aT^u$. The integral Gauss--Jacobi
identities give
\begin{align}
 \tau_a\tau_c&=J(\chi^a,\chi^c)\tau_{a+c}
          &&(a,c,a+c\not\equiv0\pmod q),\label{eq:gauss-jacobi}\\
 \tau_1\tau_{q-1}&=\ell.\label{eq:gauss-unit}
\end{align}
These follow by grouping the two summation indices according to
their sum; the zero-sum terms vanish in \eqref{eq:gauss-jacobi}
because the product character is nontrivial. They are identities
in integral cyclotomic rings and remain valid after reduction
modulo $n$; see also \cite[Theorem 2.4]{conrad}. Iteration gives
$\tau_1^q=\beta_q$.

For a shorter computation, maintain $a<q$ and $A\in R_q$ with
$\tau_1^a=A\tau_a$, starting at $a=1,A=1$. Follow the binary
expansion of $q-1$ after its leading bit. On a doubling, set
\[
 A\longleftarrow A^2J(\chi^a,\chi^a),\qquad
 a\longleftarrow2a.
\]
If the next bit is $1$, then set
\[
 A\longleftarrow A J(\chi^a,\chi),\qquad
 a\longleftarrow a+1.
\]
Every resulting index is at most $q-1$, so no step uses a
trivial product character. At the end,
$\tau_1^{q-1}=A\tau_{q-1}$; multiplying by $\tau_1$ gives
$\beta_q=\ell A$. The coefficient ring embeds in the displayed
extension, since its defining polynomial is monic, so this is an
identity in $R_q$. There are $O(\log q)$ steps and no divisions.
\end{proof}

\subsection{Combining the prime degrees}

Substitution $Z\mapsto Y^{r/q}$ defines a homomorphism
$R_q\to C_n$: the polynomial $\Phi_q(Y^{r/q})$ is divisible by
$\Phi_r(Y)$ over $\Z$. Form
\begin{equation}\label{eq:b}
 b(Y)=\prod_{q\mid r}\beta_q(Y^{r/q})^{r/q}
        \quad\text{in }C_n.
\end{equation}
Each computation in $Z$ takes place in a ring of rank $q-1$;
the substitution combines its result with the other prime degrees
in the common coefficient ring $C_n$.

\begin{proposition}[Completeness on prime inputs]\label{prop:complete}
If $n$ is prime, the construction \eqref{eq:b} satisfies all
three equalities \eqref{eq:normalize}--\eqref{eq:main}.
\end{proposition}

\begin{proof}
Put $N=n^D$. In $R_q$, $N$th powering fixes every coefficient,
because $n^D\equiv1\pmod q$ and the scalar field is $\F_n$.
For the Gauss sum from Lemma~\ref{lem:jacobi}, Frobenius and
reindexing give
\[
 \tau_1^N=\sum_{u\ne0}\chi_q(u)T^{uN}
          =\chi_q(N)^{-1}\tau_1
          =Z^{-Ds_q}\tau_1.
\]
The scalar $\ell_q$ is invertible modulo $n$, so
\eqref{eq:gauss-unit} makes $\tau_1$ a unit. Since
$\beta_q=\tau_1^q$ and $q\mid N-1$, cancellation gives
\begin{equation}\label{eq:local-normalization}
 \beta_q^{(N-1)/q}=Z^{-Ds_q}.
\end{equation}
This identity holds in $R_q$, by its injection into the monic
extension used in the proof. Substituting into \eqref{eq:b},
\[
 b^{(N-1)/r}=Y^{-D\sum_{q\mid r}s_q(r/q)}=Y.
\]
Indeed, modulo each prime $q\mid r$, all but its own summand
vanish, and \eqref{eq:sq} makes the remaining exponent $1$.
The Chinese remainder theorem gives the exponent $1$ modulo $r$.
This proves \eqref{eq:normalize}. Finally, $n$th powering is
Frobenius in every commutative algebra of characteristic the prime
$n$. It gives $b(Y)^n=b(Y^n)$ and
$(1+X)^n=1+X^n$, proving the other two equalities.
\end{proof}

\section{The algorithm and its bit complexity}\label{sec:algorithm}

Here is the complete algorithm, with $\delta=1/32$ and an effective
constant $C_\delta$ from Lemma~\ref{lem:aux} fixed in advance.
One may take $C_\delta$ to be an integer, so the search bound
\eqref{eq:search-bound} is computed by integer powering.

\begin{enumerate}
\item Reject integers below $2$; accept $2$ and $3$; reject even
      inputs greater than $2$. Reject nontrivial perfect powers.
\item Construct $D,r$ as in Proposition~\ref{prop:degrees}.
      For every prime $q\mid r$, if $q\mid n$, accept when $q=n$
      and otherwise reject. Construct $\Phi_r(Y)$ modulo $n$.
\item Sieve the primes up to $B$ in \eqref{eq:search-bound}.
      If a prime in the sieve divides $n$, accept when it equals
      $n$ and otherwise reject. For every $q\mid r$, select a prime
      $\ell_q\le B$ with $\ell_q\nmid n$, $q\mid\ell_q-1$, and
      $n^{(\ell_q-1)/q}\not\equiv1\pmod{\ell_q}$.
      Reject if any such search fails.
\item Compute the tables $j_q$, the residues $s_q$, and the
      Jacobi-sum coefficients $\beta_q$ using the binary recurrence.
      Form $b$ by \eqref{eq:b}.
\item Check \eqref{eq:normalize}, \eqref{eq:compat}, and
      \eqref{eq:main}. Accept if and only if all three hold.
\end{enumerate}

\begin{proof}[Proof of Theorem~\ref{thm:main}]
If the input is prime, preprocessing cannot reject it.
Lemma~\ref{lem:aux}, using Theorem~\ref{thm:zero-free}, ensures
all searches succeed, and Proposition~\ref{prop:complete} proves
acceptance. A composite input accepted at the last step would be
a prime power by Theorem~\ref{thm:criterion}, contrary to the
preprocessing. The earlier acceptance branches identify the input
with a known prime. This proves correctness on every input.

We give the cost bounds in some detail, including the construction
of the coefficient. The degree search costs $L^{4+o(1)}$ by
Proposition~\ref{prop:degrees}. Elementary integer-root methods
for exponents up to $L$ detect nontrivial perfect powers within
$L^{4+o(1)}$ as well.

The sieve up to $B$ uses $\softO(B)$ bit operations. Compute $n$
modulo all its primes with a product tree and a remainder tree
\cite{remainder}. The total bit length $S$ of these moduli is at
most $O(B\log B)$. With fast multiplication and division, the
trees cost $\softO(S+L)=\softO(B+L)$: at every tree level the
node lengths sum to $O(S)$, and there are $O(\log B)$ levels.
This avoids an extra factor $L$ from reducing $n$ separately for
each candidate.

Each small-prime modular exponentiation costs a polynomial in
$\log B$. There are $O(\log r)$ primes $q\mid r$, all of size
$L^{o(1)}$. Testing the candidates for all $q$ consequently costs
$L^{1/\delta+o(1)}$ in total. The character tables $j_q(u)$
are computed by modular exponentiation and lookup among $q$
possibilities. Each needed Jacobi sum is accumulated by scanning
$u=2,\ldots,\ell_q-1$ and counting its exponent modulo $q$.
These counts have $O(\log\ell_q)$ bits; the scan does not require
an $L$-bit modular addition for each term. All scans, including
the $O(\log q)$ sums for each $q$, therefore fit the same
$L^{1/\delta+o(1)}$ bound. Arithmetic in the local rings of rank
$q-1=L^{o(1)}$, after reduction of the counts, costs at most
$L^{2+o(1)}$ over all $q$.

Put $m=\varphi(r)\le r=L^{2+o(1)}$. Fast polynomial
multiplication and monic reduction give multiplication costs
\[
 \softO(mL)\quad\text{in }C_n,\qquad
 \softO(mrL)\quad\text{in }\A_n.
\]
These bounds apply over $\Z/n\Z$ without assuming it is a field.
For the second, polynomial multiplication in $X$ over $C_n$,
followed by $X^r=b$, gives the stated bound; equivalently one
may use bivariate Kronecker substitution and monic reductions.
Forming \eqref{eq:b} costs $L^{3+o(1)}$. The coefficient
checks cost $L^{4+o(1)}$, because the exponents have at most
$O(DL)$ bits. Substitution $Y\mapsto Y^n$ is performed using
$Y^r=1$ followed by reduction modulo $\Phi_r$.

The main congruence needs $O(L)$ multiplications in $\A_n$, for
total cost
\[
 \softO(mrL^2)=L^{6+o(1)}.
\]
Finally, $\Phi_r$ is constructible by exact monic divisions from
the polynomials $Y^d-1$ for the known divisors $d\mid r$.
There are $r^{o(1)}$ such divisors. Even classical polynomial
division with elementary $L$-bit coefficient arithmetic costs
at most $r^{2+o(1)}L^2=L^{6+o(1)}$. This also shows that no
polynomial-factorization subroutine is hidden in the construction.
Combining the bounds gives
\[
 L^{\max\{6,1/\delta\}+o(1)}=L^{32+o(1)},
\]
which proves the theorem.
\end{proof}

The auxiliary-prime search determines the exponent in the stated
bound. The construction does not factor $n-1$, $n^D-1$, or any
$\ell_q-1$.

\subsection{Dependence on the zero-free boundary}
\label{sec:conditional-bounds}

We now record conditional improvements under stronger analytic
hypotheses. The proof of Theorem~\ref{thm:main} is complete without
these additional inputs. For the refinements, it suffices to have
a common zero-free
boundary for the fields $\Q(\mu_{12q})$ and
$\Q(\mu_{12q},n^{1/q})$, for prime inputs $n>3$ and odd primes
$q\ne n$; a theorem over all number fields supplies this uniformly.

\begin{corollary}[Transfer of uniform zero-free bounds]
\label{cor:analytic-complexity}
Fix a real algebraic number $1/2\le\theta<1$. Suppose that the
Dedekind zeta functions of these fields have no zeros in
$\Re s>\theta$. With a corresponding choice of search cutoff,
the construction gives a deterministic primality algorithm with
bit complexity
\begin{equation}\label{eq:general-complexity}
 L^{\max\{6,\,1/(1-\theta)\}+o(1)}.
\end{equation}
In particular, a uniform zero-free boundary at $5/6$
gives $L^{6+o(1)}$. Under GRH for the stated fields, the auxiliary-prime
search and local Jacobi-sum computations cost $L^{2+o(1)}$, and the total
cost is $L^{6+o(1)}$.
\end{corollary}

\begin{proof}
For $1/2<\theta<1$, put $\delta=1-\theta$ in
Lemma~\ref{lem:aux}. Its proof uses only the two stated fields,
so the weaker hypothesis of this corollary suffices. Since all
$q\mid r$ satisfy $q\le D+1=L^{o(1)}$, the common cutoff is
$B=L^{1/(1-\theta)+o(1)}$. The sieve, remainder trees, character
tables, and local Jacobi-sum computations have this same bound; the
remaining work costs $L^{6+o(1)}$.

For completeness, the cutoff is effective even when the exponent
is not an integer. If $A$ denotes the integer base in
\eqref{eq:search-bound}, approximate $1/(1-\theta)$ from above by
a dyadic rational within $1/\log_2 A$. This increases its power
of $A$ by at most a factor of $2$. Since $\theta$ is fixed and
algebraic, such an approximation is computable in time polynomial
in $\log A$. Its denominator can be $O(\log A)$, so integer
powering and integer-root extraction compute an upper-rounded
cutoff in time polynomial in $\log A$. This cost is negligible
in the displayed bound.

At the endpoint $\theta=1/2$, Lemma~\ref{lem:explicit} remains
valid with $\delta=1/2$. The proof of Lemma~\ref{lem:aux} gives
an error $O(q(L+\log q+2)\sqrt{x})$ in
\eqref{eq:difference}, and discarded prime-power and higher-degree
terms $O(q\sqrt{x}\log^2(2x))$. Choosing a sufficiently large
absolute constant $T$ and taking
\[
 x=T\bigl(q(L+\log q+2)\bigr)^2
       \log^4\!\bigl(2q(L+\log q+2)\bigr)
\]
makes their sum, together with the ramified contribution, smaller
than the main term. Thus a suitable prime satisfies
\[
 \ell_q\ll\bigl(q(L+\log q+2)\bigr)^2
       \log^4\!\bigl(2q(L+\log q+2)\bigr).
\]
Replacing the logarithms by upper-rounded binary logarithms changes
only an absolute factor. The common cutoff is consequently
$L^{2+o(1)}$. The sieve, tables, and local coefficient arithmetic
fit this bound. Forming $b$ from the local coefficients costs
$L^{3+o(1)}$, and the remaining algebraic work costs $L^{6+o(1)}$.
GRH for these Dedekind zeta functions supplies the endpoint
hypothesis. Finally,
$1/(1-5/6)=6$ gives the stated threshold.
\end{proof}

\begin{remark}[Extent of the analytic improvement]\label{rem:analytic-limit}
Every stronger uniform zero-free boundary satisfying
Corollary~\ref{cor:analytic-complexity} yields the corresponding
bound in \eqref{eq:general-complexity} by changing the
auxiliary-prime search cutoff. The degree parameters, Jacobi-sum
construction, and algebraic soundness proof remain the same.

The exponent supplied by \eqref{eq:general-complexity} reaches $6$
at $\theta=5/6$ and stays $6$ throughout $1/2\le\theta\le5/6$.
Thus $6$ is the best exponent obtained by improving the zero-free
input alone in the present construction and cost analysis. It is
an upper-bound threshold, not a lower bound for every implementation
of this test. A smaller exponent would require a faster algebraic
stage or a different parameter construction. In particular, the
present argument does not yield a bound below the
$L^{6+o(1)}$ Gaussian-period algorithm of Lenstra and Pomerance
\cite{lenstra-pomerance}.

The RH hypothesis here concerns the Dedekind zeta functions of
the varying cyclotomic and Kummer fields. RH for $\zeta_{\Q}$
alone does not supply that hypothesis. Table~\ref{tab:analytic}
records the resulting exponents; the $o(1)$ terms tend to zero as
$L\to\infty$, with the analytic boundary fixed.
\end{remark}

\begin{table}[htbp]
\centering
\small
\begin{tabular}{@{}lccc@{}}
\toprule
Analytic input & Boundary $\theta$ & Search exponent & Total exponent \\
\midrule
Cureton--Thorner \cite{cureton-thorner} & $31/32$ & $32$ & $32$ \\
Assumed $5/6$ boundary & $5/6$ & $6$ & $6$ \\
GRH for the relevant fields & $1/2$ & $2$ & $6$ \\
\bottomrule
\end{tabular}
\caption{Exponents in the bounds $L^{a+o(1)}$. The search column
includes the character tables and local Jacobi-sum computations.
The first row is the main result; all subsequent rows are conditional
consequences of the indicated stronger analytic inputs.}
\label{tab:analytic}
\end{table}

\begin{remark}[A two-equality form]
If retaining the exponent $n$ is unnecessary, one may replace
\eqref{eq:compat} and \eqref{eq:main} by
\[
 (1+X)^{n^D}=1+YX\quad\text{in }\A_n.
\]
Together with \eqref{eq:normalize}, this gives the rotation used
in Theorem~\ref{thm:criterion} directly. Prime inputs satisfy it,
and the same proof applies. Its asymptotic verification bound is
unchanged because $D=L^{o(1)}$, although its exponent is longer.
\end{remark}

\section{A fixed-base conjecture}\label{sec:conjecture}

The coefficient $b(Y)$ in the proved construction enforces the
normalization needed for the cyclic rotation. A simpler possibility
is to fix the coefficient at $2$ and choose the degree directly
from the input. The following conjecture proposes such a test for
odd inputs, including a preliminary base-$2$ Fermat check.

Agrawal's conjecture, as recorded in \cite[Conjecture 6.1]{aks}, asserts
that if $r$ is prime, $r\nmid n$, and
\[
 (X-1)^n\equiv X^n-1\pmod{n,X^r-1},
\]
then $n$ is prime or $n^2\equiv1\pmod r$. The conjecture below
is a modification of that proposal: it uses roots of $2$ in place
of roots of unity, the shift $+1$, and the least nondivisor of
$n(n-1)$ as the degree. The two conjectures are not asserted to
be equivalent.

\begin{conjecture}[Fixed-base test]\label{conj:fixed}
Let $n>3$ be odd, and let $r$ be the least positive integer not
dividing $n(n-1)$. If
\begin{gather}
 2^{n-1}\equiv1\pmod n,\label{eq:fixed-fermat}\\
 (1+X)^n\equiv1+X^n\pmod{n,X^r-2},\label{eq:fixed-conjecture}
\end{gather}
then $n$ is prime.
\end{conjecture}

The base-$2$ Fermat congruence is an explicit condition of the test.

The degree $r$ need not be prime. It is necessarily a prime power:
otherwise it would be the least common multiple of its proper
prime-power factors, each of which divides $n(n-1)$ by minimality.

\begin{proposition}[Consequence of the fixed-base conjecture]
If Conjecture~\ref{conj:fixed} holds, its test is a deterministic
primality algorithm using $\softO(L^3)$ bit operations.
\end{proposition}

\begin{proof}
Handle inputs at most $3$ and even inputs directly, then reject
unless \eqref{eq:fixed-fermat} holds. Binary modular powering
performs this preliminary check in $\softO(L^2)$ bit operations.
By the choice of $r$,
\[
 \operatorname{lcm}(1,\ldots,r-1)\mid n(n-1).
\]
Chebyshev's estimate
$\log\operatorname{lcm}(1,\ldots,m)=\Theta(m)$ therefore gives
$r=O(\log n)=O(L)$. Testing consecutive integers for divisibility
finds $r$ within the claimed bound. Primes pass
\eqref{eq:fixed-fermat} by Fermat's little theorem and
\eqref{eq:fixed-conjecture} by Frobenius, and the conjecture
excludes composites. Multiplication in
$(\Z/n\Z)[X]/(X^r-2)$ costs $\softO(rL)$, and binary powering
uses $O(L)$ multiplications, for total $\softO(L^3)$.
\end{proof}

The open step is to obtain a comparable restriction on composite
inputs from the fixed coefficient $2$ and this degree rule. The
normalization and cyclic rotation used in Theorem~\ref{thm:criterion}
are supplied by our coefficient construction; they have not been
deduced from the congruence in Conjecture~\ref{conj:fixed}.

\begin{samepage}
\section*{AI Assistance Disclosure}

OpenAI's GPT-6 Astra developed and wrote the proofs, including the
deterministic parameter construction and complexity analysis, through
an extended mathematical dialogue with the author. The author supplied
the fixed-base test, conjectural framework, research notes, and
references; proposed applying OpenAI's recent analytic number-theoretic
results to deterministic parameter selection; and directed the
investigation of Jacobi sums, power nonresidues, composite degrees, multivariate
quotient rings, and finite-field root counting.
\end{samepage}


\begingroup\small
\begin{thebibliography}{99}

\bibitem{apr}
L. M. Adleman, C. Pomerance, and R. S. Rumely,
\emph{On distinguishing prime numbers from composite numbers},
Annals of Mathematics \textbf{117} (1983), 173--206.
Theorem 3 records the Odlyzko--Pomerance parameter construction.
\url{https://math.dartmouth.edu/~carlp/PDF/paper37.pdf}.

\bibitem{aks}
M. Agrawal, N. Kayal, and N. Saxena, \emph{PRIMES is in P},
Annals of Mathematics \textbf{160} (2004), 781--793.
\href{https://doi.org/10.4007/annals.2004.160.781}{doi:10.4007/annals.2004.160.781}.
Agrawal's conjecture is stated in Section 6; see also the authors'
version at
\url{https://www.cse.iitk.ac.in/users/manindra/algebra/primality_v6.pdf}.

\bibitem{axiom-lean}
Axiom Math, D. Cureton, and J. Thorner,
\emph{Dedekind Zeta Zero-Free Region}, Lean 4 formalization,
October 2026. The declaration \texttt{DedekindZetaZeroFree.main}
in \texttt{Solution/DedekindZetaSolution.lean} proves the
Dedekind-zeta statement; the repository reports verification
using only Lean's standard axioms.
\url{https://github.com/AxiomMath/DedekindZetaZeroFree}.

\bibitem{bernstein}
D. J. Bernstein, \emph{Proving primality in essentially quartic
random time}, Mathematics of Computation \textbf{76} (2007),
389--403.
\href{https://doi.org/10.1090/S0025-5718-06-01786-8}{doi:10.1090/S0025-5718-06-01786-8}.
Section references are to the December 3, 2004 draft:
\url{https://cr.yp.to/papers/quartic-20041203.pdf}.

\bibitem{remainder}
D. J. Bernstein, N. Heninger, and T. Lange,
\emph{FactHacks: Remainder tree}.
\url{https://facthacks.cr.yp.to/remainder.html}.

\bibitem{berrizbeitia}
P. Berrizbeitia, \emph{Sharpening ``PRIMES is in P'' for a large
family of numbers}, Mathematics of Computation \textbf{74} (2005),
2043--2059.
\href{https://doi.org/10.1090/S0025-5718-05-01727-8}{doi:10.1090/S0025-5718-05-01727-8}.

\bibitem{cheng}
Q. Cheng, \emph{Primality proving via one round in ECPP and one
iteration in AKS}, preprint, 2003.
\url{https://arxiv.org/abs/math/0301179}.

\bibitem{conrad}
K. Conrad, \emph{Gauss and Jacobi sums on finite fields and
$\Z/m\Z$}, expository notes.
\url{https://kconrad.math.uconn.edu/blurbs/gradnumthy/Gauss-Jacobi-sums.pdf}.

\bibitem{cureton-thorner}
D. Cureton and J. Thorner,
\emph{A quasi-Riemann hypothesis for Dedekind zeta functions},
Hexagon preprint 2610.00200, version 2, October 12, 2026.
Theorem 1.1 and Section 15.
\url{https://hexagonmath.org/2610.00200v2}.

\bibitem{lenstra-pomerance}
H. W. Lenstra, Jr., and C. Pomerance,
\emph{Primality testing with Gaussian periods},
Journal of the European Mathematical Society \textbf{21} (2019),
1229--1269.
\href{https://doi.org/10.4171/JEMS/861}{doi:10.4171/JEMS/861}.

\bibitem{openai-artin}
OpenAI, \emph{Primitive roots for every admissible integer base},
preprint, October 4, 2026. Theorem 1.2 and Section 9.
\url{https://github.com/openai/math/blob/main/preprints/Primitive-roots-for-every-admissible-integer-base-October-4-2026/primitive-roots-all-integer-bases.pdf}.

\end{thebibliography}
\endgroup
\end{document}